\documentclass[11pt,a4paper]{article}
\usepackage{amsmath}
\usepackage{amsfonts}
\usepackage{amssymb}
\usepackage{graphicx}
\usepackage{enumerate}
\usepackage{multicol}

\newtheorem{thm}{Theorem}[section]

\newtheorem{lemma}[thm]{Lemma}
\newtheorem{de}[thm]{Definition}

\newcommand{\QED}{\hfill$\square$}

\title {
    \bf {Distance spectra and Distance energy of \\ Integral Circulant Graphs}
}

\author
{
{\large \sc Aleksandar Ili\' c \footnotemark[3]} \\
{\em \normalsize Faculty of Sciences and Mathematics, Vi\v segradska 33, 18000 Ni\v s, Serbia} \\
{\normalsize e-mail: { \tt aleksandari@gmail.com }}
}

\begin{document}

\maketitle

\begin{abstract}
The distance energy of a graph $G$ is a recently developed energy-type invariant, defined as the
sum of absolute values of the eigenvalues of the distance matrix of $G$. There was a vast research
for the pairs and families of non-cospectral graphs having equal distance energy, and most of these
constructions were based on the join of graphs. A graph is called circulant if it is Cayley graph
on the circulant group, i.e. its adjacency matrix is circulant. A graph is called integral if all
eigenvalues of its adjacency matrix are integers. Integral circulant graphs play an important role
in modeling quantum spin networks supporting the perfect state transfer. In this paper, we
characterize the distance spectra of integral circulant graphs and prove that these graphs have
integral eigenvalues of distance matrix $D$. Furthermore, we calculate the distance spectra and
distance energy of unitary Cayley graphs. In conclusion, we present two families of pairs $(G_1,
G_2)$ of integral circulant graphs with equal distance energy -- in the first family $G_1$ is
subgraph of $G_2$, while in the second family the diameter of both graphs is three.
\end{abstract}

 \vspace{0.2cm}



\section{Introduction}

Let $G$ be a simple undirected graph with $n$ vertices. The vertices of $G$ are labeled as $0, 1,
\ldots, n - 1$. The distance between the vertices $i$ and $j$ is the length of a shortest path
between them, and is denoted by $d (i, j)$. The diameter of $G$, denoted by $diam (G)$, is the
maximum distance between any pair of vertices of $G$.

Let $A$ be the adjacency matrix of $G$, and
$\lambda_1, \lambda_2, \ldots, \lambda_n$ be the eigenvalues of the
graph $G$. The energy of $G$ is defined as the sum of absolute
values of its eigenvalues \cite{Gu78},
$$
E(G) = \sum_{i = 1}^n |\lambda_i|.
$$
The energy is a graph parameter stemming from the H\"{u}ckel molecular orbital approximation for
the total $\pi$-electron energy (for recent survey on molecular graph energy see \cite{Gu01} and
\cite{Br06}).

The distance matrix of a graph $G$ is the square matrix $D(G) = [d (i, j)]_{i,j=1}^n$. The
eigenvalues of the distance matrix $D (G)$, labeled as $\mu_1 \geqslant \mu_2 \geqslant \ldots
\geqslant \mu_n$, are said to be the distance or D-eigenvalues of $G$ and to form the distance or
$D$-spectrum of $G$ \cite{BuHa90}. The sum of distance eigenvalues is zero, $\sum_{i = 1}^n \mu_i =
0$. The characteristic polynomial and the eigenvalues of the distance matrix of a graph were
considered in \cite{GrPo71}--\cite{StIl09}.

\begin{de}
The distance energy $DE (G)$ of a graph $G$ is the sum of absolute values of the eigenvalues of the distance
matrix of G.
\end{de}

Distance energy is a useful molecular descriptor in QSPR modeling, as demonstrated by Consonni and
Todeschini in \cite{CoTo08}. To avoid trivial cases, we say that the graphs $G$ and $H$ of the same
order are D-equienergetic if $DE (G) = DE (H)$, while they have distinct spectra of distance
matrices.

Our motivation for this research came from various constructions of non-cospectral, equienergetic
graphs with equal number of vertices \cite{BrStGu04}--\cite{BoViAb08}. Indulal et al.
\cite{InGuVi08} constructed pairs of D-equienergetic graphs on $n$ vertices for $n \equiv 1 \pmod
3$ and for $n \equiv 0 \pmod 6$. Ramane et al. \cite{RaReGuWa09} proved that if $G_1$ and $G_2$ are
$r$-regular graphs on $n$ vertices and $diam (G_i) \leqslant 2$ , $i = 1, 2$, then $DE (L^k (G_1))
= DE (L^k (G_2))$ for $k \geqslant 1$, where $L^k (G)$ is the $k$-th iterated line graph of $G$. In
\cite{RaGuRe08} the authors obtain the eigenvalues of the distance matrix of the join of two graphs
whose diameter is less than or equal to two, and construct pairs of non D-cospectral,
D-equienergetic graphs on $n$ vertices for all $n \geqslant 9$. Stevanovi\' c and Indulal
\cite{StIn09} further generalized this result and described the distance spectrum and energy of the
join-based compositions of regular graphs in terms of their adjacency spectrum. Those results are
used to show that there exist a number of families of sets of non-cospectral graphs with equal
distance energy, such that for any $n \in N$, each family contains a set with at least $n$ graphs.
All these constructions from the literature are based on graph products and most of presented
graphs have diameter two.

A graph is called \emph{circulant} if it is Cayley graph on the circulant group, i.e. its adjacency
matrix is circulant. A graph is called \emph{integral} if all eigenvalues of its adjacency matrix
are integers. Integral graphs are extensively studied in the literature and there was a vast
research for specific classes of graphs with integral spectrum \cite{BaCvRaSiSt02}.

Integral circulant graphs were imposed as potential candidates for modeling quantum spin networks
with periodic dynamics. For the certain quantum spin system, the necessary condition for the
existence of perfect state transfer in qubit networks is the periodicity of the system dynamics
(see \cite{SaSeSh07}). Relevant results on this topic were given in \cite{Go08}, where it was shown
that a quantum network topology based on the regular graph with at least four distinct eigenvalues
is periodic if and only if it is integral. Various properties of integral circulant graphs were
investigated in \cite{So06}--\cite{BaIl09}.

Integral circulant graphs arise as a generalization of unitary Caley graphs,
recently studied by Klotz and Sander \cite{KlSa07}. Let $D$ be a set of positive,
proper divisors of the integer $n > 1$. Define the graph $ICG_n (D)$
to have vertex set $Z_n = \{0, 1, \ldots, n - 1 \}$ and
edge set
$$
E (ICG_n (D)) = \left \{ \{a, b \} \mid a, b \in Z_n, \ \gcd (a - b, n) \in D \right \}.
$$

In this paper our intention is to move a step towards in the investigation of the graph theoretical
properties of integral circulant graphs that are important parameters of quantum networks. Namely,
we deal with the distance matrix of integral circulant graph $ICG_n (D)$ and distance spectra. The
integral circulant graph $ICG_{10} (1)$ is shown on Figure 1, together with its distance matrix.

\begin{figure}[ht]
\begin{multicols}{2}
  \centering
  \includegraphics [width = 4.6cm]{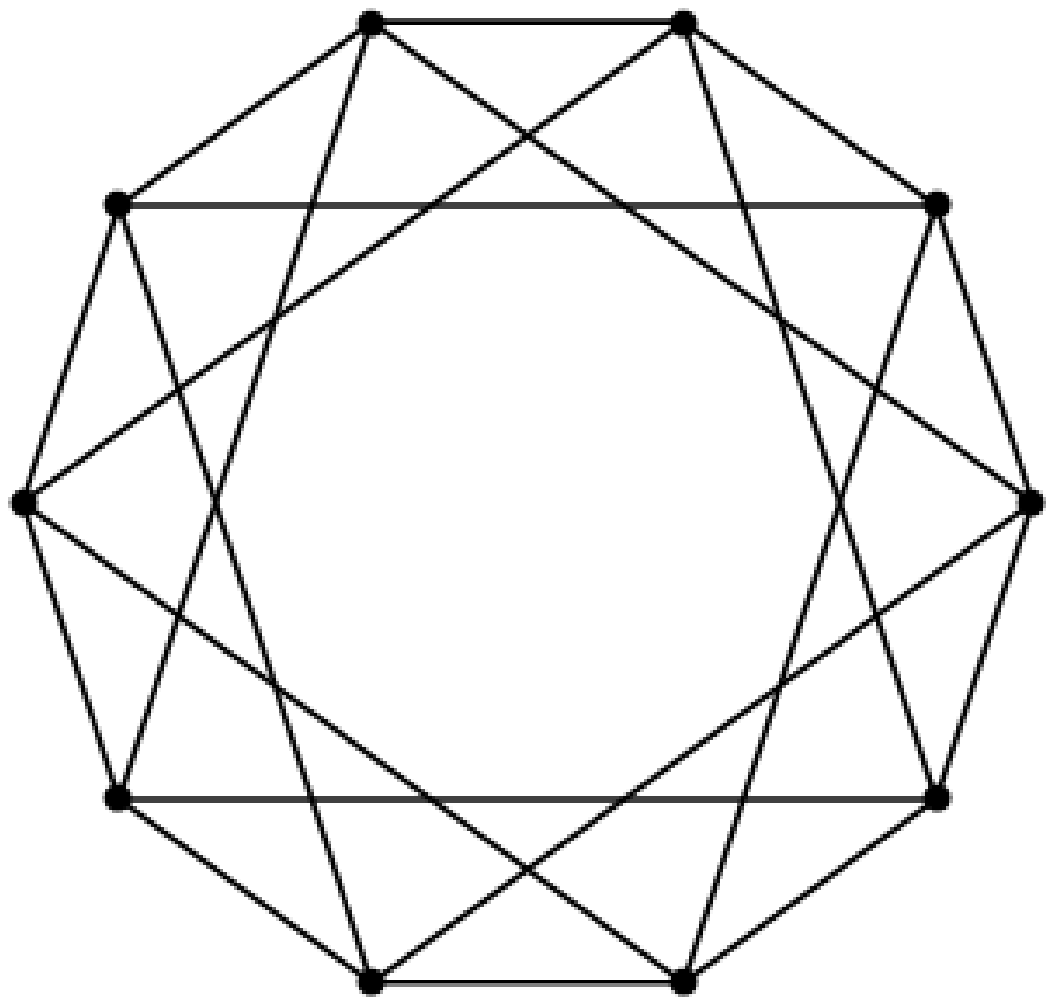}

    \columnbreak
    \centering
$$
\footnotesize D =
\begin{pmatrix}
0 &  1 &  2 &  1 &  2 &  3 &  2 &  1 &  2 &  1 \\
1 &  0 &  1 &  2 &  1 &  2 &  3 &  2 &  1 &  2 \\
2 &  1 &  0 &  1 &  2 &  1 &  2 &  3 &  2 &  1 \\
1 &  2 &  1 &  0 &  1 &  2 &  1 &  2 &  3 &  2 \\
2 &  1 &  2 &  1 &  0 &  1 &  2 &  1 &  2 &  3 \\
3 &  2 &  1 &  2 &  1 &  0 &  1 &  2 &  1 &  2 \\
2 &  3 &  2 &  1 &  2 &  1 &  0 &  1 &  2 &  1 \\
1 &  2 &  3 &  2 &  1 &  2 &  1 &  0 &  1 &  2 \\
2 &  1 &  2 &  3 &  2 &  1 &  2 &  1 &  0 &  1 \\
1 &  2 &  1 &  2 &  3 &  2 &  1 &  2 &  1 &  0 \\
\end{pmatrix}
$$
\end{multicols}
\caption { \it The integral circulant graph $ICG_{10} (1)$ and its distance matrix. }
\end{figure}

The paper is organized as follows. In Section 2 we present some preliminary results on integral
circulant graphs, while in Section 3 we prove that distance matrix of ICG graphs have integral
spectra, and characterize the vertices at distance $k$ from the starting vertex $0$. Our main
result is the description of the distance spectrum of $ICG_n (D)$. In Section 4 we calculate the
distance energy, the distance spectral radius and the Wiener index of unitary Cayley graphs. Note
that the distance energy of $ICG_n(1)$, where $n$ is even with odd prime divisor, is not fully
resolved since the sign of expression $2 (\varphi (n) - 1) - \frac{n}{2}$ can vary. Finally, in
Section 5 we present two families of non-cospectral distance equienergetic graphs $(ICG_{3p} (1),
ICG_{3p} (1, p))$ and $(ICG_{2pq} (1, p), ICG_{2pq} (1, q))$ for arbitrary primes $p > q > 3$.
These constructions are not based on the graph products, and in the first pair we have that
$ICG_{3p} (1)$ is subgraph of $ICG_{3p} (1, p)$, while in the second pair the diameter of both
graphs is three.

\section{Preliminaries}

Let us recall that for a positive integer $n$ and subset $S \subseteq \{0, 1, 2, \ldots, n - 1\}$,
the circulant graph $G(n, S)$ is the graph with $n$ vertices, labeled with integers modulo $n$,
such that each vertex $i$ is adjacent to $|S|$ other vertices $\{ i + s \pmod n \ | \ s \in S\}$.
The set $S$ is called a symbol of $G(n, S)$.

So \cite{So06} has characterized integral circulant
graphs. Let
$$
G_n (d) = \{ k\ | \ gcd(k, n)=d, \ 1 \leq k < n \}
$$
be the set of all positive integers less than $n$ having the same
greatest common divisor $d$ with $n$. Let $D_n$ be the set of
positive divisors $d$ of $n$, with $d \leqslant \frac{n}{2}$.

\begin{thm}
A circulant graph $G(n, S)$ is integral if and only if
$$
S = \bigcup_{d \in D} G_n(d)
$$
for some set of divisors $D \subseteq D_n$.
\end{thm}

It follows that the degree of $ICG_n (D)$ is equal to $\sum_{d \in D} \varphi (n / d)$, where
$\varphi (n)$ denotes the Euler phi function. Authors in \cite{BaIl09} proved that an integral
circulant graph $ICG_n (D) = ICG_n (d_1, d_2, \ldots, d_k)$ is connected if and only if $\gcd (d_1,
d_2, \ldots, d_k) = 1$.

Let $D$ be an arbitrary set of divisors $\{d_1, d_2, \ldots, d_k \}$ of $n$.
We establish the distance matrix $D_n (D)$ of integral circulant graph $ICG_n (D)$ with respect to the natural order of the
vertices $0, 1, \ldots, n - 1$. The entries $a_0, a_1, \ldots, a_{n - 1}$ of the first row of
$D_n (D)$ generate the entries of the other rows by a cyclic shift.
$$
D_n (D) = \left (
\begin{array}{cccc}
a_0 & a_1 & \ldots & a_{n - 1} \\
a_{n-1} & a_0 & \ldots & a_{n - 2} \\
\vdots & \vdots & \ddots & \vdots \\
a_1 & a_2 & \ldots & a_0 \\
\end{array}
\right )
$$

For more details on circulant matrices see \cite{Da79}. There is an explicit
formula for the eigenvalues $\mu_r$, $0 \leqslant r \leqslant n - 1$, of a circulant matrix such as $D_n$.
Define the polynomial $P_n (z)$ by the entries of the first row of $D_n$,
$$
P_n (z) = \sum_{j = 0}^{n - 1} a_j \cdot z^j
$$

The eigenvalues of $D_n$ are given by
\begin{equation}
\label{eq:ramanujan}
\mu_r = P_n (\omega^r) = \sum_{j = 0}^{n - 1} a_j \cdot \omega^{r j}, \qquad 0 \leqslant r \leqslant n - 1.
\end{equation}

Ramanujan's sum, usually denoted $c (r, n)$, is a function of two positive integer variables $r$ and $n$ defined by the formula
$$
c (r, n)= \sum_{a=1 \atop \gcd(a,n)=1}^n e^{\frac{2 \pi i}{n} \cdot a r} = \sum_{a=1 \atop \gcd(a,n)=1}^n \omega_n^{a r},
$$
where $\omega_n$ denotes a complex primitive $n$-th root of unity.
These sums take only integral values,
$$
c (r, n) = \mu \left ( \frac{n}{\gcd (r, n)} \right ) \cdot \frac{ \varphi (n) } {\varphi \left (\frac{n}{\gcd (r, n)} \right )},
$$
where $\mu$ denotes the M\" obious function. In the next section, we will use the well-known summation \cite{HaWr80}
$$
s (r, n) = \sum_{i = 0}^{n - 1} \omega_n^{i r} = \left\{
\begin{array}{l l}
  0 & \quad \mbox{ if } \quad r \nmid n\\
  n & \quad \mbox{ if } \quad r \mid n\\
\end{array} \right.
$$

In \cite{KlSa07} it was proven that \emph{gcd-graphs} (the same term as integral circulant graphs $ICG_n (D)$) have integral
spectrum,
\begin{equation}
\label{eq:eigenvalues}
\lambda_k = \sum_{d \in D} c \left(k, \frac{n}{d}\right), \qquad 0
\leqslant k \leqslant n - 1.
\end{equation}

\section{Distance spectra of ICG}

Let $N_p (i)$ be the set of vertices that are at distance $p$ from the starting vertex $i$. We will
prove that
$$
N_p (0) = \bigcup_{i = 1}^{s_p} G_n (d^{(p)}_i),
$$
for some divisor set of $n$
$$
D^{(p)} = \{d^{(p)}_1, d^{(p)}_2, \ldots, d^{(p)}_{s_p}\}.
$$

For $p = 1$, this fits the definition of integral circulant graph $ICG_n (D)$ for
divisor set $D^{(1)} = D$.
Let $a$ be an arbitrary vertex at the distance $p$ from the starting vertex $0$. Assume that
$\gcd (a, n) = d$. If follows that
there exist vertices $0 \equiv v_0, v_1, \ldots, v_{p - 1}, v_p \equiv a$, such that for all $0 \leqslant r \leqslant p - 1$,
the vertices $v_p$ and $v_{p + 1}$ are adjacent. In other words,
$$
\gcd (v_r - v_{r + 1}, n) = d_{i_r}, \qquad \mbox{for} \qquad r = 0, 1, \ldots, p - 1, \qquad 1 \leqslant i_r \leqslant k.
$$

Let $b$ be an arbitrary vertex (different than $a$) such that $\gcd (b, n) = d$. This means that there are
integers $a'$ and $b'$ relatively prime with $n$, such that $a = a' \cdot d$ and $b = b' \cdot d$.
Since $\gcd (a', n) = 1$, from B\' ezout's identity we have
$$
a' \cdot a^* \equiv 1 \pmod n,
$$
where $a^*$ is the inverse of $a'$ modulo $n$.
Now, consider the following vertices $0 \equiv u_0, u_1, \ldots, u_{p - 1}, u_p \equiv b$, defined as
$$
u_r = v_r \cdot a^* b'.
$$
For $r = 0$ we have $u_0 = 0$, and for $r = p$ we have that $u_p = a \cdot a^* b' = d \cdot
(a' \cdot a^*) \cdot b' = d \cdot b' = b$. For $r = 0, 1, \ldots, p - 1$ holds
$$
\gcd (u_r - u_{r + 1}, n) = \gcd ((v_r - v_{r + 1}) \cdot a^* b', n) = \gcd (v_r - v_{r + 1}, n) = d_{i_r}.
$$
This proves that the distance from $0$ to the vertex $b$ is less than or equal to $k$, $d (0, b)
\leqslant d (0, a)$. Similarly, we prove that $d (0, a) \leqslant d (0, b)$ and finally, we get $d
(0, a) = d (0, b) = k$.

It is easy to see that the set of all divisors is the union of all sets $D^{(p)}$,
$p = 1, 2, \ldots, diam (G)$.
Using the relation~(\ref{eq:ramanujan}), we get the distance spectra of $ICG_n (D)$:
$$
\mu_r = 1 \cdot \sum_{i = 1}^{s_1} c \left (r, \frac{n}{d^{(1)}_i} \right ) +
2 \cdot \sum_{i = 1}^{s_2} c \left (r, \frac{n}{d^{(2)}_i} \right ) + \ldots +
diam (G) \cdot \sum_{i = 1}^{s_{diam (G)}} c \left (r, \frac{n}{d^{(diam (G))}_i} \right).
$$

This proves the main result of this section.
\begin{thm}
Integral circulant graph $ICG_n (D)$, where $D$ is an arbitrary set of
divisors of $n$, has integral distance spectra.
\end{thm}

The Wiener index of $G$ is the sum of distances
between all pairs of vertices,
$$
W (G) = \sum_{a, b \in V} d (a, b).
$$
The Wiener index is considered as one of the most used topological indices with high correlation
with many physical and chemical properties of molecular compounds (for recent results and
applications of Wiener index see \cite{DoEn01}).

In \cite{In09} it is proven that $\mu_0 \geqslant \frac{2 W (G)}{n}$, with equality if
and only if all row sums are equal. Since the distance matrix of $ICG_n (D)$ is also a
circulant matrix, it follows that
$$
W (G) = \frac{n \cdot \mu_0}{2}.
$$

\section{Distance energy of UCG}

The Unitary Cayley graph $ICG_n (1) \equiv X_n$ has the vertex set $V (X_n) = Z_n$ and the edge set
$$
E (X_n) = \{ (a, b): a, b \in Z_n, \ \gcd (a - b, n) = 1 \}.
$$
The graph $X_n$ is regular of degree $\varphi (n)$. Unitary Cayley graphs are highly symmetric and
have some remarkable properties connecting graph theory and number theory. Fuchs \cite{Fu05} showed
that the maximal length of an induced cycle in $X_n$ is $2^k + 2$, where $k$ is the number of
different prime divisors of $n$. Klotz and Sander \cite{KlSa07} determined the diameter, clique
number, chromatic number and eigenvalues of unitary Cayley graphs.

In next four subsections, we will calculate the distance spectra and distance energy of $X_n$.

\subsection{$n$ is a prime number}

In this case, $X_n$ is a complete graph $K_n$ with diameter $1$.
The spectra of $X_n$ is $\{ p-1, -1, -1, \ldots, -1\}$ and
$$
DE (X_n) = 2 (n - 1).
$$

\subsection{$n$ is the power of $2$}

If $n = 2^k$ for $k > 1$, then $X_n$ is a complete bipartite graph with the vertex partition
$$
V (X_n) = \{ 0, 2, \ldots, n - 2\} \cup \{ 1, 3, \ldots, n - 1 \},
$$
and has diameter $2$. The adjacency spectrum of $X_n$ is $\{\frac{n}{2}, -\frac{n}{2}, 0, 0, \ldots, 0\}$, and consequently
the distance spectrum of $X_n$ is $\{\frac{3n}{2}-2,\frac{n}{2}-2,-2,-2,\ldots, -2\}$. Since $\varphi (2^k) = 2^{k-1}$, it follows that
$$
DE (X_n) = \left| \frac{3n}{2}-2 \right| + \left| \frac{n}{2}-2 \right| + 2 (n - 2) = 4(n - 2).
$$

\subsection{$n$ is odd composite number}

We need the following result from \cite{KlSa07},

\begin{thm}
\label{thm:neighbors}
The number of common neighbors of distinct vertices $a$ and $b$ in the unitary Cayley graph $X_n$ is given by $F_n (a - b)$,
where $F_n (s)$ is defined as
$$
F_n (s) = n \prod_{i = 1}^k \left ( 1 - \frac{\varepsilon (p)}{p} \right ), \qquad \mbox{with} \qquad
\varepsilon (p) = \left\{
\begin{array}{l l}
  1 & \quad \mbox{if} \quad p \mid s\\
  2 & \quad \mbox{if} \quad p \nmid s\\
\end{array} \right.
$$
\end{thm}

Since $2$ does not divide $n$, according to Theorem \ref{thm:neighbors}
all factors in the expansion of $F_n (a - b)$ are greater than zero.
Therefore, in this case there is a common neighbor of every pair of
distinct vertices, which implies that the diameter is $2$ and
$$
D_n (X_n) = 2 (J_n - I_n) - A_n (X_n),
$$
where $J_n$ is the matrix of ones and $I_n$ is the identity matrix. The eigenvalues of $D_n$ are given by
$$
\mu_r = 2 \sum_{i = 0}^{n - 1} \omega_n^{i r} - 2 - c (r, n).
$$
Finally, using the relation (\ref{eq:eigenvalues}) the spectra of $D_n$ is
$$
\{ 2 (n - 1) - \varphi (n), -2 - c (1, n), -2 - c (2, n), \ldots, -2 - c (n - 1, n) \}
$$
Therefore, the distance spectral radius of $X_n$ is $2 (n - 1) - \varphi (n)$.
We can calculate the distance energy of $X_n$ using similar technique as in \cite{Il09}.
The nullity of a graph $G$, denoted as $\eta (G)$, is the multiplicity of zero as the eigenvalue in adjacency spectra.

\begin{lemma}
\label{le-nullity} The nullity of $X_n$ is $n - m$, where $m = p_1
p_2 \cdot \ldots \cdot p_k$ is the maximal square-free divisor of
$n$.
\end{lemma}

Consider the following sum
$$
S = \sum_{i = 1}^n | - c (i, n) - 2 | = \sum_{i = 1}^n |c (i, n) + 2|.
$$
We already know the nullity of $X_n$, and therefore we will sum only
the non-zero eigenvalues from the spectra of $X_n$. Divide the sum
$S$ in two parts: when $\frac{n}{\gcd (n, i)}$ is a square-free
number with an even number of divisors and when $\frac{n}{\gcd (n,
i)}$ is a square-free number with an odd number of divisors. The
number of even subsets of $\{ p_1, p_2, \ldots, p_k \}$ is equal to
the number of odd subsets of $\{ p_1, p_2, \ldots, p_k \}$.

In the first case, we have
\begin{equation}
\label{eq-odd} \sum_{i \in S_1} \left ( \frac{\varphi (n)}{\varphi
(\frac{n}{\gcd (i, n)})} + 2 \right) = \varphi (n) \cdot 2^{k - 1} +
2 \sum_{l \mid m, \ \mu (l) = 1} \varphi (l).
\end{equation}
Let $l$ be a square-free number that divides $m$ with an even number
of prime factors. The number of solutions of the equation
$\frac{n}{\gcd (i, n)} = l$ is equal to $\varphi (l)$. For all $0
\leqslant i < n$ that satisfy $n = l \cdot \gcd (i, n)$ we have
$$
\frac{\varphi (n)}{\varphi (\frac{n}{\gcd (i, n)})} \cdot \varphi (l)
+ 2 \varphi (l) = \varphi(n) + 2 \varphi (l).
$$
After taking the summation for all $l \mid m$ with $\mu (l) = 1$ we
derive the identity (\ref{eq-odd}). Analogously, in the second case
we have
$$
\sum_{i \in S_2} \left ( \frac{\varphi (n)}{\varphi (\frac{n}{\gcd (i, n)})} -
2 \right ) = \varphi (n) \cdot 2^{k - 1} - 2 \sum_{l \mid m, \ \mu (l)
= -1} \varphi (l).
$$

Since Euler function $\varphi (n)$ is multiplicative, after adding
the above sums and eigenvalue $-2$ we get
\begin{eqnarray*}
S &=& 2\cdot (n - m) + \varphi (n) \cdot 2^k + 2 \sum_{l \mid m} \mu (l)
\varphi (l) \\
&=& 2n - 2m + \varphi (n) \cdot 2^k + 2 \prod_{i = 1}^k (1 - \varphi
(p_i)) = 2n - 2m + \varphi (n) \cdot 2^k + 2 \prod_{i = 1}^k (2 - p_i).
\end{eqnarray*}

Finally, the distance energy of $X_n$ equals
\begin{eqnarray*}
DE (X_n) &=& S - |\varphi (n) + 2| + |2n - \varphi (n) - 2| \\
&=& 2 \left ( 2n + \varphi (n) (2^{k - 1} - 1) - m - 2 + \prod_{i = 1}^k (2 - p_i) \right ).
\end{eqnarray*}

\subsection{$n$ is even with odd prime divisor}

As in the previous cases there are $\varphi (n)$ ones in the first row of matrix $D_n$. The vertex
$0$ is not adjacent to even vertices $2, 4, \ldots, n - 2$, but any two even vertices have a common
neighbor by Theorem \ref{thm:neighbors}. Since all neighbors of an odd vertex are even vertices, we
conclude that the number of vertices at distance $2$ is exactly $\frac{n}{2} - 1$. Similarly, any
two vertices $a$ and $b$ of $X_n$, which are both odd have a common neighbor. This proves that the
diameter of $X_n$ is $3$ and there are exactly $\frac{n}{2} - \varphi (n)$ vertices at distance $3$
from the starting vertex $0$.

It follows that the polynomial $P_n (z)$ equals
$$
P_n (z) = 3 \sum_{i = 0}^{n - 1} z^i - \sum_{i = 0}^{\frac{n}{2} - 1} z^{2i} - 2 - 2 \sum_{i = 1 \atop \gcd(i, n)=1}^n z^i.
$$

The eigenvalues of $D_n$ are given by
$$
\mu_r = 3 \sum_{i = 0}^{n - 1} \omega_n^{i r} - \sum_{i = 0}^{\frac{n}{2} - 1} \omega_n^{2 i r} - 2 - 2 \cdot c (r, n).
$$
For $r = 0$, the distance spectral radius is equal to
$$
\mu_0 = 3 n - \frac{n}{2} - 2 - 2 \varphi (n) = \frac{5n}{2} - 2 (\varphi (n) + 1).
$$
For $r = \frac{n}{2}$, it follows that $c (n/2, n) = -\varphi (n)$ and
$$
\mu_{n/2} = 0 - \frac{n}{2} - 2 - 2 \cdot c (\frac{n}{2}, n) = 2 (\varphi (n) - 1) - \frac{n}{2}.
$$
For $r \neq 0$ and $r \neq \frac{n}{2}$, the value of the sum
$$
\sum_{i = 0}^{\frac{n}{2} - 1} \omega_n^{2 i r} = \frac{\omega_n^{nr} - 1}{\omega_n^{2r} - 1}
$$
is zero. Finally, the spectra of $D_n$ consists of $\frac{5n}{2} - 2 (\varphi (n) + 1)$, $2 (\varphi (n) - 1) - \frac{n}{2}$ and
$$
-2 - 2 \cdot c (1, n), -2 - 2 \cdot c (2, n), \ldots, -2 - 2 \cdot c (n/2 - 1, n), -2 - 2 \cdot c (n/2 + 1, n), \ldots -2 - 2 \cdot c (n - 1, n).
$$

The sum $S' = \sum_{i = 0}^{n - 1} |c (i, n) + 1|$ can be computed analogously, and since the smallest prime of $n$ is $p_1 = 2$
we have
$$
S' = n - m + \varphi (n) \cdot 2^k + \prod_{i = 1}^k (2 - p_i) = n - m + \varphi (n) \cdot 2^k.
$$

The distance energy of $X_n$ is equal to
\begin{eqnarray*}
DE (X_n) &=& 2 S' - |2 + 2 \cdot c (0, n)| - |2 + 2 \cdot c (n/2, n)| + \left|\frac{5n}{2} - 2 (\varphi (n) + 1)\right| + \left|2 (\varphi (n) - 1) - \frac{n}{2}\right| \\
&=& 2n - 2m + \varphi (n) \cdot 2^{k + 1} - (2 + 2 \varphi (n)) - (2 \varphi (n) - 2) +\left(\frac{5n}{2} - 2 (\varphi (n) + 1)\right) \\
&\phantom{=}& + \left|2 (\varphi (n) - 1) - \frac{n}{2}\right|.
\end{eqnarray*}

The value of $\left|2 (\varphi (n) - 1) - \frac{n}{2}\right|$ cannot be resolved, since it takes
both positive, zero and negative values (consider for example $n = 6$, $n = 10$ and $n = 12$). By
using the representation $n = 2^k \cdot m$ with $m$ being odd, it follows $\varphi (n) = \varphi
(2^k) \varphi (m) = 2^{k-1} \varphi (m)$. Therefore, the equation $2 (\varphi (n) - 1) =
\frac{n}{2}$ is equivalent with
\begin{equation}
2 \varphi (n) = n + 1.
\end{equation}
The equation (4) is related to the still open conjecture of Lehmer \cite{Le32}, and some obvious
solutions involve prime Fermat numbers, $n = F_k - 2 = 2^{2^k} - 1$ for $k = 0, 1, 2, 3, 4, 5$.

\section{Distance equienergetic graphs}

Let $n = pq$, for odd prime numbers $p$ and $q$. According to Subsection 4.3, the distance energy
of $ICG_n (1)$ equals
\begin{eqnarray*}
DE (ICG_n (1)) &=& 3pq + 2(p -1)(q - 1) - n + 2 (2 - p)(2 - q) \\
&=& 6 (p - 1)(q - 1).
\end{eqnarray*}

Consider the graph $ICG_n (1, p)$ -- its diameter is two, and only the vertices satisfying $\gcd (a
- b, n) = q$ are at distance $2$. Using the formula (\ref{eq:eigenvalues}) we get the distance
eigenvalues of $ICG_n (1, p)$,
$$
\mu_r = c (r, pq) + c (r, q) + 2 \cdot c (r, p).
$$
Based on the greatest common divisor of $r$ and $pq$, it follows that
$$
\mu_r = \left\{
\begin{array}{l l}
  1 - 1 + 2 \cdot (-1) & \quad \mbox{ if } p \nmid r \ \mbox{ and } \ q \nmid r\\
  - \varphi(p) - 1 + 2 \cdot \varphi (p) & \quad \mbox{ if } p \mid r \ \mbox{ and } \ q \nmid r\\
  - \varphi (q) + \varphi (q) + 2 \cdot (-1) & \quad \mbox{ if } p \nmid r \ \mbox{ and } \ q \mid r\\
  \varphi (pq) + \varphi (q) + 2 \cdot \varphi (p) & \quad \mbox{ if } p \mid r \ \mbox{ and } \ q \mid r\\
\end{array} \right.
$$
Finally, the spectra of $ICG_n (1, p)$ consists of $-2$ with multiplicity $pq - q$, $p - 2$ with multiplicity $q - 1$,
and $pq + p - 2$ with multiplicity $1$. Therefore, the distance energy equals
\begin{eqnarray*}
DE (ICG_n (1, p)) &=& q (p - 1) \cdot 2 + (q - 1) \cdot (p - 2) + 1 \cdot (pq + p - 2) \\
&=& 4 q (p - 1).
\end{eqnarray*}

Only for $q = 3$ and $n = 3p$, we have the identity $DE (ICG_n (1)) = DE (ICG_n (1, p))$, which
represents a family of distance equienergetic integral circulant graphs. Note that $ICG_{3p} (1)$
is a subgraph of $ICG_{3p} (1, p)$ for all primes $p > 3$.

Let now $n = 2pq$, where $p$ and $q$ are different odd prime numbers.
We will prove that graphs $ICG_n (1, p)$ and $ICG_n (1, q)$ have equal distance energy,
$$
DE (ICG_n (1, p)) = DE (ICG_n (1, q)) = 4 (3pq - p - q - 1).
$$

According to Subsection 4.4, the distance eigenvalues of bipartite $ICG_n (1, p)$ are given by
$$
\mu_r = c (r, 2pq) + c (r, 2q) + 2 c(r, pq) + 2 c (r, p) + 2 c (r, q) + 3 c (r, 2p) + 3 c (r, 2).
$$
Based on the greatest common divisor of $r$ and $2pq$, it follows that
$$
\mu_r = \left\{
\begin{array}{l l}
  -1 + 1 + 2 - 2 - 2 + 3 - 3 & \quad \mbox{ if } 2 \nmid r \ \mbox{ and } \ p \nmid r \ \mbox{ and } \ q \nmid r\\
  1 - 1 + 2 - 2 - 2 - 3 + 3 & \quad \mbox{ if } 2 \mid r \ \mbox{ and } \ p \nmid r \ \mbox{ and } \ q \nmid r\\
  (p - 1) + 1 - 2 (p - 1) + 2 (p - 1) - 2 - 3 (p - 1) - 3 & \quad \mbox{ if } 2 \nmid r \ \mbox{ and } \ p \mid r \ \mbox{ and } \ q \nmid r\\
  -(p - 1) - 1 - 2 (p - 1) + 2 (p - 1) - 2 + 3 (p - 1) + 3 & \quad \mbox{ if } 2 \mid r \ \mbox{ and } \ p \mid r \ \mbox{ and } \ q \nmid r\\
  (q - 1) - (q - 1) - 2 (q - 1) - 2 + 2 (q - 1) + 3 - 3 & \quad \mbox{ if } 2 \nmid r \ \mbox{ and } \ p \nmid r \ \mbox{ and } \ q \mid r\\
  -(q - 1) + (q - 1) - 2 (q - 1) - 2 + 2 (q - 1) - 3 + 3 & \quad \mbox{ if } 2 \mid r \ \mbox{ and } \ p \nmid r \ \mbox{ and } \ q \mid r\\
  pq - 2p - 2 & \quad \mbox{ if } 2 \nmid r \ \mbox{ and } \ p \mid r \ \mbox{ and } \ q \mid r\\
  3pq + 2p - 2 & \quad \mbox{ if } 2 \mid r \ \mbox{ and } \ p \mid r \ \mbox{ and } \ q \mid r\\
\end{array} \right.
$$

$$
\mu_r = \left\{
\begin{array}{l l}
  3pq + 2p - 2 & \quad \mbox{ for } r = 0\\
  pq - 2p - 2 & \quad \mbox{ for } r = pq\\
  -2p-2 & \quad \mbox{ for }\gcd (r, 2q) = 1\\
  2p-2 & \quad \mbox{ for }\gcd (r, q) = 1\\
  -2 & \quad \mbox{ otherwise }
\end{array} \right.
$$
Since the multiplicity of distance eigenvalues $2p-2$ and $-2p-2$ is $q - 1$, we finally get
\begin{eqnarray*}
DE (ICG_n (1, p)) &=& (3pq + 2p - 2) + (pq - 2p - 2) + (q - 1)(2p + 2) + (q - 1)(2p - 2) + 2(pq - 2q) \\
&=& 12pq - 4p - 4q - 4,
\end{eqnarray*}
which is a symmetric expression involving $p$ and $q$. Note that the diameter of graphs $ICG_{2pq}
(1, p)$ and $ICG_{2pq} (1, q)$ is equal to three.

These constructions of infinite families of distance equienergetic graphs are the first ones
derived not using the product of graphs nor iterated line graphs.
It would be of great interest to calculate the distance energy of other classes
of integral circulant graphs, and explore new families of non-cospectral distance
equienergetic integral graphs.

\vspace{0.5cm} {\bf Acknowledgement. } This work was supported by the research grant 144007 of the
Serbian Ministry of Science and Technological Development. The author is grateful to the anonymous
referees for their comments and suggestions.

\end{document}